\documentclass[11pt]{article}

\title{Analysis of a Double Kruskal Theorem}
\author{Timothy Carlson}
\date{}

\usepackage{amsmath}    % need for subequations
\usepackage{amssymb}
\usepackage{graphicx}   % need for figures
\usepackage{verbatim}   % useful for program listings
\usepackage{color}      % use if color is used in text
\usepackage{subfigure}  % use for side-by-side figures
\usepackage{hyperref}   % use for hypertext links, including those to external documents and URLs                                    

\newtheorem{prop}{Proposition}[section]
\newtheorem{cor}[prop]{Corollary}
\newtheorem{thm}[prop]{Theorem}
\newtheorem{lem}[prop]{Lemma}

\newtheorem{dfn}[prop]{Definition}

\newcommand{\qed}{\hfill {\bf QED}}
\newcommand{\blankline}{\vspace{4 mm}}
\newcommand{\halfblankline}{\vspace{2 mm}}

\newcommand{\bfP}{\mbox{${\bf P}$}}

\newcommand{\bfQ}{\mbox{${\bf Q}$}}

\newcommand{\sfT}{\mbox{${\sf T}$}}
\newcommand{\sfS}{\mbox{${\sf S}$}}

\newcommand{\preq}{\mbox{$\, \preceq\, $}}
\newcommand{\preqcov}{\mbox{$\, \preceq^c \, $}}

\newcommand{\leqone}{\mbox{$\leq_1$}}
\newcommand{\leqtwo}{\mbox{$\leq_2$}}

\begin{document}
\maketitle

\begin{center}
  The Ohio State University, Columbus, OH 43210 USA \\
  email: carlson@math.ohio-state.edu
\end{center}

\blankline

%: abstract

\noindent
{\small 
{\bf Abstract.}
The strength of an extension of Kruskal's Theorem [\ref{rkr}] to certain pairs of cohabitating trees is calibrated showing that it is independent of the theory $\Pi^1_1-{\bf CA}_0$ or, equivalently, ${\bf KP}\ell_0$.
}

\blankline

This paper is a sequel to Carlson [\ref{rca}] where we considered whether families of cohabitating trees are wqo under inf preserving embeddings. 
We found that when considering families of three trees, the natural candidates are not wqo. A result of Laver [\ref{rla76}] places strong restrictions on families of pairs of cohabitating trees which are wqo by showing that the family of pairs of cohabitating linear orderings is not wqo.

The following concepts arise naturally in the author's work in proof theory.

A structure
$$(X,\leqone,\leqtwo)$$
is a {\bf double forest}
if
\begin{itemize}
\item
Both $(X, \leqone)$ and $(X,\leqtwo)$ are finite forests (i.e. finite partial orderings in which the set of predecessors of any element is linearly ordered).
\item
For all $a,b\in X$, 
$$a\leqtwo b \ \ \ \Longrightarrow \ \ \ a\leqone b$$
\item
For all $a,b,c\in X$,
$$ a\leqone b \leqone c \ {\it and } \ a\leqtwo c \ \ \ \Longrightarrow \ \ \ a\leqtwo b$$
\end{itemize}
The {\bf height} of a nonempty double forest $(X,\leq_1,\leq_2)$ is the height of the forest $(X,\leq_2)$ i.e. the natural number $n$ such that the size of the longest chain is $n+1$.
If we strengthen the first condition in the definition of double forest by requiring that $(X,\leq_1)$ and $(X,\leq_2)$ are trees (i.e. forests with a minimum element, called the {\bf root}), then we call $(X,\leqone,\leqtwo)$ a {\bf double tree}.

Assume that $(Q,\preq)$ is a quasiordering i.e. reflexive and transitive. Recall that $(Q,\preq)$ is a {\bf well quasiordering}  if $(Q,\preq)$ has no bad sequences (where an infinite sequence $q_0,q_1,\ldots,q_n,\ldots$ is {\bf bad} iff it is not the case that there are $i<j$ such that $q_i\preq q_j$).

The two theorems below follow from the stronger results from Carlson [\ref{rca}] where the collection of double trees is replaced by the collection of pure patterns of order 2.

\blankline

\noindent {\bf Theorem.} Assume $n$ is a natural number. The collection of double trees of height at most $n$ is wqo under inf preserving embeddings. 

\blankline

While the full collection of double trees is not wqo under embeddings, they are wqo under a weaker notion.

 Let $(X,\leqone,\leqtwo)$ and $(X^*,\leq^*_1,\leq^*_2)$ be double forests. 
An injection $h$ of $X$ into $X^*$ is a {\bf covering} if for $i=1,2$
$$a\leq_i b \ \ \ \Longrightarrow \ \ \ h(a)\leq_i^* h(b)$$
for all $a,b\in X$.

\blankline

\noindent {\bf Double Kruskal Theorem.} The collection of finite double trees is  wqo under coverings.  

\blankline

We will consider the strength of several variants of the theorems above.
In particular, our investigation will show that the Double Kruskal Theorem is equivalent over ${\bf ACA}_0$ to the uniform $\Pi^1_1$ reflection principle for $\Pi^1_1-{\bf CA}_0$ and is independent over $\Pi^1_1-{\bf CA}_0$ (or, equivalently, ${\bf KP}\ell_0$).
Our approach is similar to that taken in Simpson [\ref{rbu}].
In particular, we will make use of a system of ordinal notations due to W. Buchholz.

While we will reference the formal theories ${\bf RCA}_0$, ${\bf WKL}_0$, ${\bf ACA}_0$,  ${\bf ID}_n$, $\Pi^1_1-{\bf CA}_0$  and ${\bf KP}\ell_0$ along with various notions from proof theory, the first three sections which include the key results, Theorem \ref{tmain} and Corollary \ref{cmain}, can be read without any specialized knowledge.  For the applications in the last two sections, we assume enough familiarity with ${\bf RCA}_0$ to recognize that the proofs in earlier sections can be carried out in ${\bf RCA}_0$ and that, except where stated otherwise, the arguments in the last two sections can also be carried out in ${\bf RCA}_0$.
Each case were we go beyond ${\bf RCA}_0$, a result from the literature will be referenced.
 The reader will only need to know that the theories ${\bf RCA}_0$, ${\bf WKL}_0$, ${\bf ACA}_0$ and $\Pi^1_1-{\bf CA}_0$ are successively increasing in strength and that ${\bf ID}_n$ implies the first order part of ${\bf ACA}_0$.
Simpson [\ref{rsi10}] is our reference for ${\bf RCA}_0$.

\section{Buchholz Notations}
\label{sbu}

In [\ref{rbu}], Buchholz defines the $\psi$-functions and appeals to them to define a system of ordinal notations. 
We will review those notations in this section.

Fix a sequence $D_0,D_1,\ldots,D_n,\ldots,D_\omega$ of formal symbols. 
We will also treat $0$ as a formal symbol in the following definition.

\begin{dfn}
Define  a set $T$ of formal terms and a subset $P$ of $T$ inductively by the following clauses.
\begin{itemize}
\item[{\rm (T1)}]
$0\in T$.
\item[{\rm (T2)}]
If $a_0,\ldots,a_n\in P$ where $n\geq 1$ then $(a_0,\ldots,a_n)\in T$.
\item[{\rm (T3)}]
If $u\leq \omega$ and $a\in T$ then $D_ua\in P$ (and, hence, $D_ua\in T$).
\end{itemize}
The elements of $P$ are called {\bf principal terms}.
The {\bf order} of a principal term $D_ua$ is $u$.
For $a=(a_0,\ldots,a_n)\in T$ with $1\leq n$, define the order of $a$ to be the order of $a_0$.
Also define the order of 0 to be 0.
For $a\in T$, write $ord(a)$ for the order of $a$.
\end{dfn}

The parentheses and commas are to be treated as formal symbols in the definition of $T$ and $P$ though we will also use $(a_0,\ldots,a_n)$ to represent the sequence with components $a_0,\ldots,a_n$ at times.

In [\ref{rbu}],  the elements of $T$ are interpreted as ordinals using the $\psi$-functions. We will not need this interpretation here, but we mention that the symbol $0$ represents the ordinal $0$, $(a_0,\ldots,a_n)$  represents the ordinal which is the sum of the ordinals represented by the $a_i$ and $D_ua$ represents the application of $\psi_u$ to the ordinal represented by $a$.

Notice that by focusing on the subscripts of the symbols $D_u$ which occur in an element $a$ of $T$, $a$ can be seen as a finite planar forest whose nodes are tagged by ordinals $u$ with $u\leq \omega$ ($0$ can be seen as the empty forest).
The following ordering on $T$ can then be seen as the natural generalization of the lexicographic ordering on sequences of ordinals bounded by $\omega$.

\begin{dfn}
The binary relation $\prec$ on $T$ is defined inductively by the following clauses:
\begin{itemize}
\item[{\rm ($\prec$1)}] 
For $b\in T$, $0\prec b$ iff $b\not=0$.
\item[{\rm ($\prec$2)}] 
For $a,b\in T$ and $u,v\leq\omega$,  $D_ua\prec D_v b$ iff either  $u<v$ or both $u=v$ and $a\prec b$.
\item[{\rm ($\prec$3)}]
Assume $a,b\in T$.
\begin{itemize}
\item[(a)]
If $a\in P$ and $b=(b_0,\ldots,b_n)$ where $1\leq n$ then $a\prec b$ iff $a\preq b_0$ (i.e. $a\prec b_0$ or $a=b_0$).
\item[(b)]
If $a=(a_0,\ldots,a_m)$ where $1\leq m$ and $b\in P$ then $a\prec b$ iff $a_0\prec b$.
\item[(c)]
 If  $a=(a_0,\ldots,a_m)$ and $b=(b_0,\ldots,b_n)$ where $1\leq m,n$ then $a\prec b$ iff one of the following holds
\begin{itemize}
\item[{\rm (i)}]
$m<n$ and $a_i=b_i$ for $i\leq m$.
\item[{\rm (ii)}]
There exists $i\leq n,m$ such that $a_j=b_j$ for $j<n$ and $a_i\prec b_i$.
\end{itemize}
\end{itemize}
\end{itemize}
\end{dfn}

We will sometimes identify elements of $T$ with sequences of elements of $P$ as follows:
\begin{itemize}
\item
$0$ is identified with the empty sequence $()$.
\item
For $u\leq\omega$ and $a\in T$, $D_ua$ is identified with the sequence $(D_ua)$.
\item
For $1\leq n <\omega$ and $a_0,\cdots,a_n\in P$, the formal term $(a_0,\ldots,a_n)$ is identified with the sequence $(a_0,\ldots,a_n)$.
\end{itemize}
Under this identification, the ordering on $T$ agrees with the lexicographic ordering derived from the ordering above restricted to $P$.

The relation $\prec$ is easily seen to be a strict linear ordering of $T$.
Moreover, if $a,b\in T$ and $ord(a)<ord(b)$ then $a\prec b$.

The interpretation of elements of $T$ as ordinals mentioned above is not an injection.
Therefore,  $\prec$ is not the ordering inherited from the ordinals.
In fact, one easily sees that $\prec$ is not a well-ordering.
These facts will not be used later.

We will restrict $\prec$ to a subset $OT$ of $T$ below.

\begin{dfn}
For $u\leq \omega$, $G_ua$ is defined for $a\in T$ inductively so that
\begin{itemize}
\item[{\rm (G1)}]
$G_u0=\emptyset$.
\item[{\rm (G2)}]
If $n\geq 1$ and $a_0,\ldots,a_n\in P$ then 
$$G_u(a_0,\ldots,a_n) = G_ua_0\cup\cdots \cup G_ua_n$$
\item[{\rm (G3)}]
If $v\leq\omega$ and $a\in T$ then
\begin{equation*}
G_uD_va \ = \ 
\begin{cases}
\emptyset & \text{if $v<u$}, \\
\{a\}\cup G_ua & \text{if $u\leq v$}.
\end{cases}
\end{equation*}
\end{itemize}
\end{dfn}

\begin{dfn}
The set $OT$ of {\bf ordinal terms} is defined inductively by
\begin{itemize}
\item[{\rm (OT1)}]
$0\in OT$.
\item[{\rm (OT2)}]
If $n\geq 1$, $a_0,\ldots,a_n\in OT$ are principal terms and $a_n\preq \cdots \preq a_0$ then $(a_0,\ldots,a_n)\in OT$.
\item[{\rm (OT3)}]
If $u\leq \omega$ and $a\in OT$ with $G_ua\prec a$ then $D_ua\in OT$.
\end{itemize}
For $u\leq \omega$, $OT(u)$ is the collection of ordinal terms $a$ such that $v<u$ whenever $D_v$ occurs in $a$.
\end{dfn}

\blankline 

Clearly, if $a\in OT$ and $u\leq\omega$ then $G_ua\subseteq OT$.

Notice that if $(a_0,\ldots,a_n)\in OT$ then $ord(a_0)\geq\cdots\geq ord(a_n)$.

\section{Monotone Double Forests and Collapsing}
\label{smo}

In this section, we consider double forests whose nodes are tagged by natural numbers so as to descend along $\leq_1$.
We will also define a preliminary version of collapsing functions, analogues of the $\psi$-functions, on certain finite sequences of such tagged double forests.

\begin{dfn}
A {\bf monotone double forest} is a pair $({\bf P},\delta)$ where ${\bf P}$ is a double forest and $\delta$ maps $|{\bf P}|$, the underlying set for $\bf P$, into $\omega$ such that
$$x\leq_1^{\bf P} y \ \ \ \Longrightarrow \ \ \ \delta(x)\geq \delta(y)$$
and
$$x\leq_2^{\bf P}y \ \ \ \Longrightarrow \ \ \ \delta(x)=\delta(y)$$
Define $M2F$ to be the collection of monotone double forests.
Let $M2T$ be the collection of elements $({\bf P},\delta)$ of $M2F$
such that  $(|{\bf P}|,\leq_1^{\bf P})$ is a tree.
For $x\in |{\bf P}|$, the {\bf order} of $x$ in $({\bf P}, \delta)$  is $\delta(x)$.
For $({\bf P},\delta)\in M2T$, the {\bf root} of $({\bf P},\delta)$ is the minimal element of $\leq_1^{\bf P}$ and the {\bf order} of $({\bf P},\delta)$, $ord(({\bf P},\delta))$, is the order of the root.
\end{dfn}

Notice that any chain in $\leq_2^{\bf P}$ must consist of elements all of which have the same order.

We define notions like {\it embedding} and {\it isomorphism} between elements of $M2F$ as usual.

When ${\sf T}=({\bf P},\delta)\in M2F$, we will often write $|{\sf T}|$, $\leq_1^{\sf T}$, $\leq_2^{\sf T}$ and $\delta^{\sf T}$ for $|{\bf P}|$, $\leq_1^{\bf P}$, $\leq_2^{\bf P}$ and $\delta$ respectively.

\begin{dfn}
Assume $\sfT_i=({\bf P}_i,\delta_i)\in M2F$ for $i=1,2$.
 We define a function $h$ from $|\sfT_1|$ into $|\sfT_2|$ to be a {\bf covering} of $\sfT_1$ into $\sfT_2$ if $h$ is a covering of ${\bf P}_1$ into ${\bf P}_2$ and 
$$\delta_1(x)\leq \delta_2(h(x))$$
for all $x\in |\sfT_1|$.
We define $\sfT_1\, \preq^c \ \sfT_2$ if there is a covering of $\sfT_1$ into $\sfT_2$.
\end{dfn}

Notice that when $\sfT=({\bf P},\delta)\in M2T$ has order $0$, $\delta$ must be identically $0$ and we can identify \sfT\ with $\bf P$. 
In this way, the collection of elements of $M2T$ of order $0$ is identified with the collection of double forests ${\bf P}$ where $\leq_1^{\bf P}$ is a tree.
Moreover, our two notions of covering coincide under this identification.

We will eventually define an element $\sfT(a)$ of $M2T$ for each principal $a\in OT(\omega)$.
For this, we will want to define {\it collapsing} operations which interpret the symbols $D_u$ for $u<\omega$.
The first step of this definition is to define a basic collapsing operation which adds a new root below a sequence of elements  of $M2T$ of the same positive order and reduces the order by 1.

\begin{dfn}
\label{dcoll}
Assume $u\in \omega$ and ${\sf T}_1,\ldots,{\sf T}_n\in M2T$  have order $u+1$. 
If $\sfT_1,\ldots,\sfT_n$ are pairwise disjoint, define $coll({\sf T}_1,\ldots,{\sf T}_n)$ to be a structure $((X,\leq_1,\leq_2),\delta)$ which satisfies  the following conditions.
\begin{enumerate}
\item
$X=|{\sf T}_1|\cup\cdots\cup |{\sf T}_n|\cup \{r\}$ where $r\not\in |\sfT_1|\cup\cdots\cup|\sfT_n|$.
\item
Assume $1\leq i\leq n$ and $x\in |{\sf T}_i|$.
For $k=1,2$ and all $y\in X$, $x\leq_k y$ iff $x\leq_k^{{\sf T}_i}y$.
\item
For all $y\in X$,  $r\leqone y$.
\item
For all $y\in X$, $r\leq_2 y$ iff either $y=r$ or there exists $i$ such that $1\leq i\leq n$, $y\in |{\sf T}_i|$ and $\delta^{{\sf T}_i}(y)=u+1$.
\item
For all $i$ with $1\leq i\leq n$ and all $x\in |{\sf T}_i|$, $\delta(x)=min\{\delta^{{\sf T}_i}(x),u\}$.
\item
$\delta(r)=u$.
\end{enumerate}
If $\sfT_1,\ldots,\sfT_n$ are not pairwise disjoint, define 
$$coll(\sfT_1,\ldots,\sfT_n)=coll(\sfS_1,\ldots,\sfS_n)$$
where  $\sfS_1,\ldots,\sfS_n$ are chosen to be pairwise disjoint,   $\sfS_i\cong \sfT_i$ for $i=1,\ldots,n$ and $\sfS_1=\sfT_1$.
\end{dfn}

Since we will mainly be concerned with elements of $M2F$ up to isomorphism, the choice of $r$ and the choice of $\sfS_1,\ldots,\sfS_n$ in the above definition are not important. 
The condition $\sfS_1= \sfT_1$ is a technical convenience.
It could be dropped without changing the isomorphism type of the resulting structure.

\begin{dfn}
For $\sfT\in M2F$ and $x\in |\sfT|$, define $\sfT^x$ to be the substructure of $\sfT$ whose universe is the collection of $y\in|\sfT|$ such that $x\leq_1^{\sf T} y$.
\end{dfn}

\begin{lem}
\label{lcoll}
Assume $u<\omega$ and $\sfT_1,\ldots,\sfT_m\in M2T$ are pairwise disjoint and have order $u+1$.
\begin{enumerate}
\item
$coll(\sfT_1,\ldots,\sfT_n)$ is an element of $M2T$ of order $u$.
\item
Assume $1\leq i \leq m$.
If $coll(\sfT_1,\ldots,\sfT_n)=({\bf P},\delta)$ and $\sfT_i=({\bf P}_i,\delta_i)$ then ${\bf P}_i$ is a substructure of $\bf P$.
\item
If $1\leq i \leq m$ and $x\in |\sfT_i|$ has order at most $u$ in $\sfT_i$  then 
$$coll(\sfT_1,\ldots,\sfT_m)^x=\sfT_i^x$$
%\item
%If $\sfS_i\cong \sfT_i$ for $i=1,\ldots,m$ then
%$$coll(\sfS_1,\ldots,\sfS_m)\cong coll(\sfT_1,\ldots,\sfT_m)$$
\end{enumerate}
\end{lem}
{\bf Proof.}
Straightforward noticing that if ${\sf T}\in M2T$ has order $u+1$ then $\delta^{\sf T}(x)\leq u+1$ for all $x\in |{\sf T}|$ and the collection of $x$ with $\delta^{\sf T}(x)=u+1$ is closed downward in $\leq_1^{\sf T}$.
\qed

\begin{dfn}
Assume $\sfT_1,\ldots,\sfT_m\in M2F$ are pairwise disjoint. Define  $\oplus(\sfT_1,\ldots,\sfT_m)$ to be the usual disjoint union $((X,\leqone,\leqtwo),\delta)$ described by the following conditions. 
\begin{enumerate}
\item
$X=|\sfT_1|\cup\cdots\cup|\sfT_m|$
\item
$\leqone=\leq_1^{{\sf T}_1}\cup\cdots\cup\leq_1^{{\sf T}_m}$
\item
$\leqtwo=\leq_2^{{\sf T}_1}\cup\cdots\cup\leq_2^{{\sf T}_m}$
\item
$\delta=\delta^{{\sf T}_1}\cup\cdots\cup \delta^{{\sf T}_m}$
\end{enumerate}
\end{dfn}

\begin{lem}
If $\sfT_1,\ldots,\sfT_m\in M2F$ are pairwise disjoint then $\oplus(\sfT_1,\ldots,\sfT_m)\in M2F$.
\end{lem}
{\bf Proof.} 
Clear.
\qed

\begin{lem}
\label{lcases}
Assume $u<\omega$, $\sfT_1,\ldots,\sfT_m,\sfS_1,\ldots,\sfS_n\in M2T$ have order $u+1$ and both sequences $\sfT_1,\ldots,\sfT_m$ and $\sfS_1,\ldots,\sfS_n$ are pairwise disjoint.
If $coll(\sfT_1,\ldots,\sfT_m)\preq^c \ coll(\sfS_1,\ldots,\sfS_n)$ then either 
$$\oplus(\sfT_1,\ldots,\sfT_m)\ \preq^c \ \oplus(\sfS_1,\ldots,\sfS_n)$$
or there exist $j$ with $1\leq j\leq n$ and $x\in |\sfS_j|$ of order $u$ such that 
$$coll(\sfT_1,\ldots,\sfT_m)\preq^c \ \sfS_j^x$$
\end{lem}
{\bf Proof.}
Write ${\sf A}_1$ for $coll(\sfT_1,\ldots,\sfT_m)$ and ${\sf A}_2$ for $coll(\sfS_1,\ldots,\sfS_m)$.
Let $r_i$ be the root of ${\sf A}_i$ for $i=1,2$.
Let $h$ be a covering of ${\sf A}_1$ into ${\sf A}_2$.

\halfblankline

{\bf Case 1.} 
$r_2\leq_2^{{\sf A}_2}h(r_1)$.

We claim that a covering of $\oplus(\sfT_1,\ldots,\sfT_m)$ into $\oplus(\sf S_1,\ldots,\sfS_n)$ is obtained by restricting $h$. 

Assume $1\leq i \leq m$ and $x\in |\sfT_i|$.
Let $j$ be such that $h(x)\in |\sfS_j|$.
We must show $\delta^{{\sf T}_i}(x)\leq \delta^{{\sf S}_j}(h(x))$.
Since $h$ is a covering of ${\sf A}_1$ into ${\sf A}_2$, $\delta^{{\sf A}_1}(x)\leq \delta^{{\sf A}_2}(h(x))$.
By the definition of $coll$, this is equivalent to $min\{\delta^{{\sf T}_i}(x),u\}\leq min\{\delta^{{\sf S}_j}(h(x)),u\}$.
If $\delta^{{\sf T}_i}(x)\leq u$, this immediately implies the desired inequality.
So, we may assume $\delta^{{\sf T}_i}(x)=u+1$.
By definition of ${\sf A}_1$, $r_1 \leq_2^{{\sf A}_1}x$.
Since $h$ is a covering, $h(r_1) \leq_2^{{\sf A}_2}h(x)$.
Since we are assuming $r_2\leq_2^{{\sf A}_2}h(r_1)$ in this case, $r_2\leq_2^{{\sf A}_2}h(x)$.
By definition of ${\sf A}_2$, $\delta^{{\sf S}_j}(h(x))=u+1$.

The rest of the proof that the restriction of $h$ is a covering follows from part 2 of Lemma \ref{lcoll}.

\halfblankline

{\bf Case 2.} 
$r_2 \not\leq_2^{{\sf A}_2}h(r_1)$.

We will show $x=h(r_1)$ witnesses the second disjunct in the conclusion of the lemma.

Since $r_2\not\leq_2^{{\sf A}_2}h(r_1)$, there exists $j$ such that $h(r_1)\in |\sfS_j|$ and, moreover, $\delta^{{\sf S}_j}(h(r_1))\leq u$.
Since $h$ is a covering and $r_1$ has order $u$ in ${\sf A}_1$, $h(r_1)$ has order at least $u$ in ${\sf A}_2$.
Equivalently, $u\leq min\{\delta^{{\sf S}_j}(h(r_1)), u\}$.
Therefore, $u\leq \delta^{{\sf S}_j}(h(r_1))$ establishing that $\delta^{{\sf S}_j}(h(r_1))=u$.

Since $h$ is a covering of ${\sf A}_1$ into ${\sf A}_2$, $h$ is a covering of ${\sf A}_1$ into ${\sf A}_2^{h(r_1)}$. 
Since $h(r_1)$ is an element of $\sfS_j$ of order at most $u$, ${\sf A}_2^{h(r_1)}=\sfS_j^{h(r_1)}$ by part 3 of Lemma \ref{lcoll}.
\qed

\blankline

For $u<\omega$, we want to define the collapsing operation $\Psi_u$ on any sequence $\sfT_1,\ldots,\sfT_n\in M2T$ with $ord(\sfT_1)\geq\cdots\geq ord(\sfT_n)$. We will refer to such sequences as being {\it order descending}. 
The following definition deals with the case when $ord(\sfT_i)>u$ for all $i$.

\begin{dfn}
\label{dpsi1}
Assume $u\in \omega$ and ${\sf T}_1,\ldots,\sfT_n\in M2T$ is order descending with 
$ord(\sfT_n)> u$.
Inductively on $ord(\sfT_1)-u$ define
$$\Psi_u(\sfT_1,\ldots,\sfT_n)\ = \ coll(\Psi_{u+1}(\sfT_1,\ldots,\sfT_i),\sfT_{i+1},\ldots,\sfT_n)$$
where $i$ is maximal such that either $i=0$ or $ord(\sfT_i)>u+1$.
 \end{dfn}
 
 Notice that $ord(\sfT_{i+1})=\cdots=ord(\sfT_n)=u+1$ in the notation of the definition.
 
 In the case where $i=0$, the interpretation of the definition is
 $$\Psi_u(\sfT_1,\ldots,\sfT_n)=coll(\sfT_1,\ldots,\sfT_n)$$
 In the case where $i=n$, the interpretation of the definition is
 $$\Psi_u(\sfT_1,\ldots,\sfT_n)=coll(\Psi_{u+1}(\sfT_1,\ldots,\sfT_n))$$
 
 % An easy induction shows that if $u$ and $\sfT_1,\ldots,\sfT_n$ are as in the definition and $\sfS_i\cong \sfT_i$ for $i=1,\ldots,n$ then $\Psi_u(\sfT_1,\ldots,\sfT_n)\cong \Psi_u(\sfS_1,\ldots,\sfS_n)$.
 
\begin{lem}
If $u<\omega$ and $\sfT_1,\ldots,\sfT_n\in M2T$ is order descending with $ord(\sfT_n)>u$ then  $\Psi_u(\sfT_1,\ldots,\sfT_n)$ is an element of $ M2T$  of order $u$.
\end{lem}
{\bf Proof.}
Simple induction on $ord(\sfT_1)-u$ using part 1 of Lemma \ref{lcoll}.
\qed

\blankline

While we will not use the following lemma directly, it is helpful in understanding  $\Psi_u(\sfT_1,\ldots,\sfT_n)$ in case $ord(\sfT_n)>u$.

\begin{lem}
\label{lPsipicture}
Assume $u,v<\omega$, $\sfT_1,\ldots,\sfT_n\in M2T$ is order descending, 
$v=ord(\sfT_1)$, $ord(\sfT_n)>u$
and
$$\Psi_u(\sfT_1,\ldots,\sfT_n)=((X,\leqone,\leqtwo),\delta)$$
There exist disjoint $\sfS_1,\ldots,\sfS_n$ such that $\sfS_i\cong\sfT_i$  for $i=1,\ldots,n$ and distinct $r_{u+1},r_{u+2},\ldots,r_v\in X-(|\sfS_1|\cup\cdots\cup|\sfS_n|)$  such that
\begin{enumerate}
\item
$X=|\sfS_1|\cup\cdots\cup |\sfS_n|\cup \{r_{u+1},r_{u+2}\ldots,r_v\}$
\item
Assume $1\leq i \leq n$ and $x\in |\sfS_i|$.
For all $y\in X$, $x\leqone y$ iff $x\leq_1^{{\sf S}_i}y$.
\item
Assume $u< w \leq v$.
For all $y\in X$,  $r_w\leqone y$ iff either $y=r_t$ for some $t$ with $w\leq t \leq v$ or $y\in |\sfS_i|$ for some $i$ such that $ord(\sfS_i)\geq w$.
\item
Assume  $1\leq i \leq n$ and $x\in |\sfS_i|$.
For all $y\in X$, $x\leq_2 y$ iff $x\leq_2^{{\sf S}_i}y$.
\item
Assume $u<w\leq v$.
For all $y\in X$, $r_w\leq_2 y$ iff either $y=r_t$ for some $t$ with $w\leq t \leq v$ or $y\in |\sfS_i|$ for some $i$ with $\delta^{{\sf S}_i}(y)\geq w$.
\item
Whenever $1\leq i \leq n$ and  $x\in |\sfS_i|$, $\delta(x)=min\{\delta^{{\sf S}_i}(x),u\}$.
\item
For $i=u+1,u+2,\ldots,v$, $\delta(r_i)=u$.
\item
$\Psi_u({\sf T})\in M2T$ has order $u$.
\end{enumerate}
\end{lem}
{\bf Proof.}
Tedious but straightforward induction on $v-u$.
\qed

\blankline

In the statement of the lemma, $r_v,\ldots,r_{u+2},r_{u+1}$ are the new roots obtained by successivly applying $coll$.  
When visualizing $\Psi_u(\sfT_1,\ldots,\sfT_n)$ with respect to \leqone\  one might view the chain 
$$r_{u+1}<_2 r_{u+2} <_2\cdots <_2 r_v $$ 
(which is closed downward with respect to \leqone) as the spine with  $|\sfS_i|$ branching away from the spine at $r_w$ where $w$ is the order of $\sfS_i$.

We will need one more operation before completing the definition of $\Psi_u$.
It will also add a new root with a specified label below a collection of  elements of $M2T$.
Notice that the first three conditions are the same as those in Definition \ref{dcoll}.

\begin{dfn}
\label{dexp}
Assume $\sfT_i \in M2T$ for $i=1,\ldots,n$ and $u<\omega$. 
If $\sfT_1,\ldots,\sfT_n$ are pairwise disjoint,  define $exp_u(\sfT_1,\ldots,\sfT_n)$ to be a structure $((X,\leq_1,\leq_2),\delta)$ which satisfies the following conditions.
\begin{enumerate}
\item
$X=|\sfT_1|\cup \cdots \cup |\sfT_n|\cup \{r\}$ where $r\not\in |\sfT_i|$ for $i=0,\ldots,n$.
\item
Assume $1\leq i\leq n$ and $x\in |\sfT_i|$.
For $k=1,2$ and all  $y\in X$, $x\leq_k y$ iff $x\leq_k^{{\sf T}_i}y$.
\item
For all $y\in X$,  $r\leqone y$.
\item
For all $y\in X$, $r\leq_2 y$ iff $y=r$.
\item
Assume $1\leq i\leq n$.
For all $x\in |{\sf T}_i|$, $\delta(x)=\delta^{{\sf T}_i}(x)$.
\item
$\delta(r)=u$.
\end{enumerate}
When $\sfT_1,\ldots,\sfT_n$ are not pairwise disjoint, define $exp_u(\sfT_1,\ldots,\sfT_n)$ to be $exp_u(\sfS_1,\ldots,\sfS_n)$ for some pairwise disjoint $\sfS_1,\ldots,\sfS_n \in M2T$ such that $\sfS_i\cong \sfT_i$ for $i=1,\ldots,n$.
\end{dfn}

\begin{lem}
\label{lexp}
Assume $u<\omega$ and $\sfT_1,\ldots,\sfT_n\in M2T$  are pairwise disjoint  of order at most $u$.
\begin{enumerate}
\item
 $exp_u({\sf T}_1,\ldots,\sfT_n)$ is an element of $ M2T$ of order $u$.
 \item For $i=1,\ldots,n$, ${\sf T}_i$ is a substructure of $exp_u(\sfT_1,\ldots,\sfT_n)$.
 \item
 If $1\leq i \leq m$ and $x\in |\sfT_i|$ then  $exp_u(\sfT_1,\ldots,\sfT_m)^x=\sfT_i^x$.
 \end{enumerate}
\end{lem}
{\bf Proof.}
Straightforward.
\qed

\blankline

We are now ready to complete the definition of $\Psi_u$.

\begin{dfn}
\label{dpsi2}
Assume $u<\omega$, $\sfT_1,\ldots,\sfT_n\in M2T$ is order descending and $ord(\sfT_n)\leq u$.
Define
$$\Psi_u(\sfT_1,\ldots,\sfT_n)\ = \ exp_u(\Psi_u(\sfT_1,\ldots,\sfT_i),\sfT_{i+1},\ldots,\sfT_n)$$
where $i$ is maximal such that either $i=0$ or $ord(\sfT_i)>u$. 
\end{dfn}

In case $i=0$, the interpretation of the definition is
$$\Psi_u(\sfT_1,\ldots,\sfT_n)=exp_u(\sfT_1,\ldots,\sfT_n)$$

\section{Reducing Ordinal Terms to Double Trees}
\label{sre}

In this section, we will define $\sfT(a)\in M2T$ for each principal  $a\in OT(\omega)$ so that
$$\sfT(a)\,  \preqcov \,\sfT(b) \ \ \ \Longrightarrow \ \ \ a\preq b$$
This will allow us to reduce arbitrary elements of $OT(\omega)$ to double trees in an analagous way.

We first need to make some observations about the notations from Section \ref{sbu}.

\begin{lem}
If $D_ua\in OT$ and $u\leq v$ then $D_va\in OT$.
\end{lem}
{\bf Proof.}
Assume $D_ua\in OT$ and $u\leq v$. Clearly, $G_vb\subseteq G_ub$ for any $b\in T$. 
Since $G_va\subseteq G_ua\prec a$, we see $D_va\in OT$.
\qed

\begin{lem}
Assume $D_u(a_1,\ldots,a_n)\in OT$.
If $1\leq m<n$ then $D_u(a_1,\ldots,a_m)\in OT$.
\end{lem}
{\bf Proof.}
Arguing by contradiction, there exists $(b_1,\ldots,b_k)\in G_u(a_1,\ldots,a_m)$ with $(a_1,\ldots,a_m)\preq (b_1,\ldots,b_k)$. 
Since $G_u(a_1\ldots,a_m)=G_ua_1\cup \cdots \cup G_ua_m\subseteq G_ua_1\cup \cdots \cup G_ua_n=G_u(a_1,\ldots,a_n)$, we have $(b_1,\ldots,b_k)\in G_u(a_1,\ldots,a_n)$.
Since $D_u(a_1,\ldots,a_n)\in OT$, this implies $(b_1,\ldots,b_k)\prec (a_1,\ldots,a_n)$.
Since $(a_1,\ldots,a_m)\preq (b_1,\ldots,b_k) \prec (a_1,\ldots,a_n)$, we conclude that $m\leq k$ and $a_i=b_i$ for $1\leq i\leq k$.
Therefore, $(a_1,\ldots,a_m,b_{m+1},\ldots,b_k)$ is a subterm of $(a_1,\ldots,a_m)$ -- contradiction.
\qed

\blankline

The following definition is from [\ref{rbu}].
Informally,  the collection of $u$-subterms of $a\in T$ is the collection of subterms of $a$  which are not in the scope of some $D_v$ for $v<u$.

\begin{dfn}
Assume $u\leq\omega$.
For $a\in T$, define the collection of {\bf $u$-subterms} of $a$ inductively by the following clauses.
\begin{enumerate}
\item
The only $u$-subterm of $0$ is $0$.
\item
Assume $a=D_vb$ where $b\in T$ and $v<u$.
The only $u$-subterm of $a$ is $a$.
\item
Assume $a=D_vb$ where $b\in T$ and $u\leq v$.
For $c\in T$, $c$ is a $u$-subterm of $a$ iff either  $c=a$ or $c$ is a $u$-subterm of $b$.
\item
Assume $a=(b_0,\ldots,b_n)$ where $n\geq 1$.
For $c\in T$, $c$ is a $u$-subterm of $a$ iff either $c=a$ or $c$ is a $u$-subterm of $b_i$ for some $i\leq n$.
\end{enumerate}
\end{dfn}

Notice that for $a\in T$, the elements of $G_ua$ are those $b\in T$ such that $D_vb$ is a $u$-subterm of $a$ for some $v\geq u$.

\begin{lem}
\label{lusubterm}
Assume $D_ua\in OT$. If $b$ is a $u$-subterm of $D_ua$ of order $u$ then $b\preq D_ua$.
\end{lem}
{\bf Proof.}
Assume $b$ is a $u$-subterm of $D_ua$ of order $u$.
If $b=D_ua$ we are done, so we may assume $b\not=D_ua$.
This assumption implies $b$ is a $u$-subterm of $a$.
There exists $c\in OT$ such that $b=D_u c$.
By the comment preceding the lemma, $c\in G_ua$.
Since $c\in G_ua \prec a$, $D_uc\prec D_ua$.
\qed

\blankline

For the following definition, we identify $D_ua$ with $D_u(a)$ when $a$ is a principal ordinal term.

\begin{dfn}
Define $\sfT(a)\in M2T$ for principal $a\in OT(\omega)$ inductively so that $|\sfT(D_u0)|=\{r\}$ for some $r$ where $r$ has order $u$ and 
$$\sfT(D_u(a_1,\ldots,a_n))= \Psi_u(\sfT(a_1),\ldots,\sfT(a_n))$$
for $D_u(a_1,\ldots,a_n)\in OT$ with $(a_1,\ldots,a_n)\not=0$.
\end{dfn}

\begin{lem}
 $ord(\sfT(a))=ord(a)$ for principal $a\in OT(\omega)$.
 \end{lem}
 {\bf Proof.}
 Immediate from the definition of $\sfT(a)$.
 \qed

\begin{lem}
\label{lT(a)}
Assume that $a=D_u(a_1,\ldots,a_m)\in OT(\omega)$.
\begin{enumerate}
\item
If $ord(a_m)\leq u$  then
$$\sfT(a)=exp_u(\sfT(D_u(a_1,\ldots,a_i)),\sfT(a_{i+1}),\ldots,\sfT(a_m))$$
where $i$ is maximal such that either $i=0$ or $ord(a_i)>u$.
\item
If $ord(a_m)>u$ then
$$\sfT(a)=coll(\sfT(D_{u+1}(a_1,\ldots,a_i)),\sfT(a_{i+1}),\ldots,\sfT(a_m))$$
where $i$ is maximal such that either $i=0$ or $ord(a_i)>u+1$.
\end{enumerate}
\end{lem}
{\bf Proof.}
The lemma follows immediately from the previous definition and Definitions \ref{dpsi1} and \ref{dpsi2} with similar interpretations of the equations e.g. the equation in part 2 is to be interpreted as 
$$\sfT(a)=coll(\sfT(a_1),\ldots,\sfT(a_m))$$
in case $i=0$.
\qed

\begin{lem}
\label{lminsubterm}
Assume $a\in OT(\omega)$ is principal. Suppose $v< \omega$ and let $X$ be the collection of $x\in |\sfT(a)|$ of order $v$.
If $x$ is a minimal element of $X$ with respect to $\leq_1^{{\sf T}(a)}$ then there exists a principal $v$-subterm $b$ of $a$  of order $v$ such that $\sfT(a)^x\cong \sfT(b)$.
\end{lem}
{\bf Proof.}
We will argue by induction on the cardinality of $|\sfT(a)|$ for principal $a\in OT(\omega)$. 

Suppose $a\in OT(\omega)$ is principal and the lemma holds with $a$ replaced by $b$ whenever $b\in OT(\omega)$ is principal and $card(\sfT(b))<card(\sfT(a))$.
Also, suppose $v<\omega$, $X$ is the collection of elements of $|\sfT(a)|$ of order $v$ and $x$ is a minimal element of $X$ with respect to $\leq_1^{{\sf T}(a)}$. 

Let $u$ be the order of $a$. 
Since $x$ has order $v$ in $\sfT(a)$, we have $v\leq u$.
If $v=u$ then $x$ is the root of $\sfT(a)$ and we can take $b=a$ in the conclusion of the lemma.
Therefore, we may assume $v<u$.

The assumption that $v<u$ implies that $a\not=D_u0$.
Therefore, $a=D_u(a_1,\ldots,a_m)$ for some $(a_1,\ldots,a_m)\in OT(\omega)$.

\halfblankline

{\bf Case 1.}
$ord(a_m)\leq u$.

Let $i$ be maximal such that either $i=0$ or $ord(a_i)>u$.
By part 1 of Lemma \ref{lT(a)},
$$\sfT(a)=exp_u(\sfT(D_u(a_1,\ldots,a_i)),\sfT_{i+1},\ldots,\sfT_m)$$
for some $\sfT_{i+1},\ldots,\sfT_m\in M2T$ such that $\sfT(D_u(a_1,\ldots,a_i)),\sfT_{i+1},\ldots,\sfT_m$ are pairwise disjoint and $\sfT_k\cong \sfT(a_k)$ for $i+1\leq k \leq m$.
Since $v<u$, either $x\in |\sfT(D_u(a_1,\ldots,a_i))|$ or $x\in |\sfT_k|$ for some $k$ with $i+1\leq k \leq m$.

We first consider the subcase where $x\in |\sfT(D_u(a_1,\ldots,a_i))|$. 
By the induction hypothesis, there is a $v$-subterm $b$ of $D_u(a_1,\ldots,a_i)$ of order $v$ such that $\sfT(D_u(a_1,\ldots,a_i))^x\cong \sfT(b)$.
By part 3 of Lemma \ref{lexp}, $\sfT(D_u(a_1,\ldots,a_i))^x=\sfT(a)^x$ so that $\sfT(a)^x\cong \sfT(b)$.
Since $v<u$, $b$ is a $v$-subterm of $a_j$ for some $j$ with $1\leq j \leq i$.
This implies that $b$ is a $v$-subterm of $a$.

The proof of  the subcase where $x\in |\sfT(a_k)|$ for some $k$ with $i+1\leq k \leq m$ is similar.

\halfblankline

{\bf Case 2.}
$ord(a_m)>u$.

Let $i$ be maximal such that either $i=0$ or $ord(a_i)>u+1$.
By part 2 of Lemma \ref{lT(a)},
$$\sfT(a)=coll(\sfT(D_{u+1}(a_1,\ldots,a_i)),\sfT_{i+1},\ldots,\sfT_m)$$
for some $\sfT_{i+1},\ldots,\sfT_m\in M2T$ such that $\sfT(D_{u+1}(a_1,\ldots,a_i)),\sfT_{i+1},\ldots,\sfT_m$ are pairwise disjoint and $\sfT_k\cong \sfT(a_k)$ for $i+1\leq k \leq m$.
Since $v<u$, either $x\in |\sfT(D_{u+1}(a_1,\ldots,a_i))|$ or $x\in |\sfT_k|$ for some $k$ with $i+1\leq k \leq m$.

The rest of the proof is analgous to Case 1 using part 3 of Lemma \ref{lcoll} instead of part 3 of Lemma \ref{lexp}.
\qed

\blankline

Assume $(a_1,\ldots,a_m),(b_1,\ldots,b_n)\in OT$ and $p: \{1,\ldots,m\}\rightarrow \{1,\ldots,n\}$.
In the proof of the following theorem, we will use the following observation: 
\begin{itemize}
\item[] 
If  for $i=1,\ldots,m$ 
$$a_i\preq b_{p(i)}$$
and whenever there exists $j\not=i$ such that $p(j)=p(i)$ 
 $$a_i\not= b_{p(i)}$$
    then 
$$(a_1,\ldots,a_m)\preq (b_1,\ldots,b_n)$$
\end{itemize}

Of course, this applies generally to descending sequences in lexicographical orderings.
Moreover, if $p$ is not injective then the conclusion can be strengthened to $(a_1,\ldots,a_m)\prec (b_1,\ldots,b_n)$.

\begin{thm}
\label{tmain}
Assume $a,b\in OT(\omega)$ are principal. If $\sfT(a) \preqcov  \sfT(b)$ then $a\, \preq \, b$.
\end{thm}
{\bf Proof.} We will argue by induction on the cardinality of $|\sfT(b)|$ with a subsidiary induction on the cardinality of $|\sfT(a)|$.

Assume $a,b\in OT(\omega)$ are principal and the conclusion of the lemma holds when replacing $a$ and $b$ by principal $a',b'\in OT(\omega)$ whenever either $card(|\sfT(b')|)<card(|\sfT(b)|)$ or both $card(|\sfT(b')|)=card(|\sfT(b)|)$ and $card(|\sfT(a')|)<card(|\sfT(a)|)$.

Assume $\sfT(a)\, \preq^c \ \sfT(b)$.
We must show $a\preq b$.

Since $\sfT(a)\, \preq^c \ \sfT(b)$, we have $ord(a)=ord(\sfT(a))\leq ord(\sfT(b))=ord(b)$.
If $ord(a)<ord(b)$ then $a\prec b$ and we are done. 
So, we may assume $ord(a)=ord(b)$.
Let $u$ be the common value.

Since $\sfT(D_u0)$ consists of a single element of order $u$, it follows easily from $\sfT(a)\, \preq^c \ \sfT(b)$  that $a\preq b$  if $a$ or $b$ is $D_u0$ (notice also that $card(\sfT(D_uc))>1$ whenever $c\in OT(\omega)$ and $c\not=0$).
Therefore, we may assume $a=D_u(a_1,\ldots,a_m)$ and $b=D_u(b_1,\ldots,b_n)$.
It suffices to show  $(a_1,\ldots,a_m)\preq(b_1,\ldots,b_n)$ in order to establish $a\preq b$.

We will consider four cases depending on whether $u<ord(a_m)$ and whether $u<ord(b_n)$.

\blankline

{\bf Case 1.} $u<ord(a_m),ord(b_n)$.

Let $i$ be maximal such that either $i=0$ or $a_i>u+1$ and let $j$ be maximal such that either $j=0$ or $b_j>u+1$.
By part 2 of Lemma \ref{lT(a)}, 
$$\sfT(a)=coll(\sfT(D_{u+1}(a_1,\ldots,a_i)),\sfT_{i+1},\ldots,\sfT_m)$$
 for some $\sfT_{i+1},\ldots,\sfT_m\in M2T$ such that $\sfT_k\cong \sfT(a_k)$ for $i+1\leq k \leq m$ and $ \sfT(D_{u+1}(a_1,\ldots,a_i)),\sfT_{i+1},\ldots,\sfT_m$ are pairwise disjoint.
 Also, 
$$\sfT(b)=coll(\sfT(D_{u+1}(b_1,\ldots,b_j)),\sfS_{j+1},\ldots,\sfS_n)$$
 for some $\sfS_{j+1},\ldots,\sfS_n\in M2T$ such that $\sfS_k\cong \sfT(b_k)$ for $j+1\leq k \leq n$ and $\sfT(D_{u+1}(b_1,\ldots,b_n)),\sfS_{j+1},\ldots,\sfS_n$ are pairwise disjoint.
 
 We will define ${\sf A}$ to be $\sfT(D_{u+1}(a_1,\ldots,a_i))$  in case $i\not=0$ and define ${\sf B}$ to be $\sfT(D_{u+1}(b_1,\ldots,b_j))$ in case $j\not=0$. 
 
 By Lemma \ref{lcases}, the three subcases which follow are exhaustive.
 
 \blankline
 
 {\bf Subcase 1.1.}
  There exists $x\in |{\sf B}|$ (so $j\not=0$) of order $u$ such that $\sfT(a)\,\preq^c \ {\sf B}^x$.
 
 Without loss of generality, $x$ is minimal among the elements of $\sf B$ of order $u$ with respect to $\leq_1^{\sf B}$.
 By Lemma \ref{lminsubterm}, ${\sf B}^x\cong \sfT(c)$ for some $u$-subterm $c$ of $D_{u+1}(b_1,\ldots,b_j)$ of order $u$.
 By the induction hypothesis, $a\preq c$. 
 By definition, $c$ must be a $u$-subterm of $(b_1,\ldots,b_j)$ implying it is a $u$-subterm of $b=D_u(b_1,\ldots,b_n)$. 
 By Lemma \ref{lusubterm}, $c\preq b$ (in fact, $c\prec b$) implying $a\preq b$.
  
 \blankline
 
 {\bf Subcase 1.2.} 
 There exists $k$ with $j+1\leq k \leq n$ and $x\in |\sfT(b_k)|$ of order $u$ such that $\sfT(a)\, \preq^c \ \sfT(b_k)^x$. 
 
Similar to the argument for Subcase 1.1, $a\preq c$ for some $u$-subterm $c$ of $b_k$ of order $u$. 
Since $b_k$ is a $u$-subterm of $b$, $c$ is a $u$-subterm of $b$.
By Lemma \ref{lusubterm}, $c\preq b$ implying $a\preq b$.

\blankline

{\bf Subcase 1.3.} $\oplus({\sf A},\sfT_{i+1},\ldots,\sfT_m)\, \preq^c \ \oplus({\sf B},\sfS_{j+1},\ldots,\sfS_n)$.

We first consider the case when $i=0$ so that $\sf A$ is undefined and $ord(a_1)=\cdots=ord(a_m)=u+1$.
If $j>0$ then $ord(b_1)>u+1$ implying $a_1\prec b_1$ and, hence, $(a_1,\ldots,a_m)\prec (b_1,\ldots,b_n)$. 
So, we may assume $j=0$. 
Let $h$ be a covering which witnesses the assumption of the subcase.
Notice that for each $k$ with $1\leq k \leq m$ there exists $p(k)$ such that $h$ maps $|\sfT_k|$ into $|\sfS_{p(k)}|$.
By the induction hypothesis, $a_k\preq b_{p(k)}$.
Moreover, if there exists $l\not=k$ with $p(l)=p(k)$ then $h$ cannot map $\sfT_k$ onto $\sfS_{p(k)}$ (since $\sfT_l$ is also mapped into $\sfS_{p(k)}$) implying $a_k\prec b_{p(k)}$ (since $\sfT_k\cong \sfT(a_k)$ has smaller cardinality than $\sfS_{p(k)}\cong \sfT(b_{p(k)})$, we have $\sfT(a_k)\not=\sfT(b_{p(k)})$ which implies $a_k\not=b_{p(k)}$).
By the observation preceding the theorem, $(a_1,\ldots,a_m)\preq (b_1,\ldots,b_n)$.

We now consider the case when $i>0$ and $\sf A$ is defined.

If $(a_1,\ldots,a_i)\prec (b_1,\ldots,b_j)$ then, since $a_{i+1}$ has order $u+1$ if it exists, $(a_1,\ldots,a_m)\prec (b_1,\ldots,b_n)$.
So, we may assume $(b_1,\ldots,b_j)\preq (a_1,\ldots,a_i)$. 

We begin by showing $(a_1,\ldots,a_i)\preq(b_1,\ldots,b_j)$ which implies $(a_1,\ldots,a_i)=(b_1,\ldots,b_j)$.

The assumption of this subcase implies there is a covering of $\sf A$ into one of ${\sf B},\sfS_{j+1},\ldots,\sfS_n$. 
This implies $D_{u+1}(a_1,\ldots,a_i)\preq D_{u+1}(b_1,\ldots,b_j)$ or $D_{u+1}(a_1,\ldots,a_i)\preq b_k$ for some $k$ with $j+1\leq k \leq n$ by the induction hypothesis. 
Our immediate goal is equivalent to $D_{u+1}(a_1,\ldots,a_i)\preq D_{u+1}(b_1,\ldots,b_j)$, so assume  $D_{u+1}(a_1,\ldots,a_i)\preq b_k$ where $j+1\leq k \leq n$.
Since $b_k$ has order $u+1$, it has the form $D_{u+1}c$ for some $c\in OT(\omega)$.
Since $D_{u+1}(a_1,\ldots,a_i)\preq D_{u+1}c$, $(a_1,\ldots,a_i)\preq c$.
Clearly, $c\in G_ub$ implying $c\prec (b_1,\ldots,b_n)$. 
Therefore, $(a_1,\ldots,a_i)\prec (b_1,\ldots,b_n)$ which implies $(a_1,\ldots,a_i)\preq (b_1,\ldots,b_j)$ (since $b_{j+1}$ has order $u+1$ if it exists).

Since $(a_1,\ldots,a_i)=(b_1,\ldots,b_j)$ which implies ${\sf A}={\sf B}$, a covering witnessing the hypothesis of this subcase can easily be modified to map $|\sf A|$ onto $|\sf B|$ thus witnessing $\oplus(\sfT_{i+1},\ldots,\sfT_m)\, \preq^c \ \oplus(\sfT_{j+1},\ldots,\sfT_n)$.
By an argument similar to that in the first paragraph of this subcase, we see $(a_{i+1},\ldots,a_m)\preq (b_{j+1},\ldots,b_n)$.
Combined with $(a_1,\ldots,a_i)=(b_1,\ldots,b_j)$, we conclude $(a_1,\ldots,a_m)\preq (b_1,\ldots,b_n)$. 

\blankline

{\bf Case 2.} $ord(a_m)\leq u < ord(b_n)$.

If $ord(a_1)\leq u$ then $a_1\prec b_1$ implying $(a_1,\ldots,a_m)\prec (b_1,\ldots,b_n)$.
So, we may assume $ord(a_1)>u$.
Let $i$ be maximal such that $ord(a_i)>u$. 
We have $\sfT(a)=exp_u(\sfT(D_u(a_1,\ldots,a_i)), \sfT(a_{i+1}),\ldots,\sfT(a_m)))$.
Since $\sfT(D_u(a_1,\ldots,a_i))$ is a substructure of $\sfT(a)$, the restriction of a covering of $\sfT(a)$ into $\sfT(b)$ is a covering of $\sfT(D_u(a_1,\ldots,a_i))$ into $\sfT(b)$ which is not onto. 
The induction hypothesis implies that $D_u(a_1,\ldots,a_i)\prec b=D_u(b_1,\ldots,b_n)$.
Therefore, $(a_1,\ldots,a_i)\prec (b_1,\ldots,b_n)$.
Since $ord(a_{i+1})\leq u$ if $i+1\leq m$, $(a_1,\ldots,a_m)\prec (b_1,\ldots,b_n)$.

\blankline

{\bf Case 3.}
$ord(b_n)\leq u < ord(a_m)$.

Let $j$ be maximal such that $j=0$ or $ord(b_j)>u$.
By part 1 of Lemma \ref{lT(a)}, $\sfT(b)=exp_u(\sfT(D_u(b_1,\ldots,b_j)),\sfT(b_{j+1}),\ldots,\sfT(b_n))$.

Since $ord(a_m)>u$, $\sfT(a)$ is obtained by an application of $coll$. 
Therefore, letting $r_1$ be the root of $\sfT(a)$, there is $x\in \sfT(a)$ such that $r\leq_2^{{\sf T}(a)}x$ and $r\not=x$.
Since $\sfT(b)$ is obtained by an applicaton of $exp_u$, letting $r_2$ be the root of $\sfT(b)$, there is no $x\in \sfT(b)$ such that $r_2\leq_2^{{\sf T}(b)}x$ and $r_2\not = x$.
Therefore,  no covering of $\sfT(a)$ into $\sfT(b)$ maps $r_1$ to $r_2$. 
Since $\sfT(a)\, \preq^c \ \sfT(b)$, either $\sfT(a)\, \preq^c \ \sfT(D_u(b_1,\ldots,b_i))$ or $\sfT(a)\, \preq^c \ \sfT(b_k)$ for some $k$ with $j+1\leq k \leq n$.
In the former case, the induction hypothesis implies $a\preq D_u(b_1,\ldots,b_i)\prec b$.
In the latter case, the induction hypothesis implies $a\preq b_k\prec b$ (since $b_k$ is a $u$-subterm of $b$ of order at most $u$).

\blankline

{\bf Case 4.}
$ord(a_m),ord(b_n)\leq u$.

Let $i$ be maximal such that either $i=0$ or $a_i>u$ and let $j$ be maximal such that either $j=0$ or $b_j>u$.
By part 1 of Lemma \ref{lT(a)}, 
$$\sfT(a)=exp_u(\sfT(D_u(a_1,\ldots,a_i)),\sfT_{i+1},\ldots,\sfT_m)$$
 for some $\sfT_{i+1},\ldots,\sfT_m\in M2T$ such that $\sfT_k\cong \sfT(a_k)$ for $i+1\leq k \leq m$ and $ \sfT(D_u(a_1,\ldots,a_i)),\sfT_{i+1},\ldots,\sfT_m$ are pairwise disjoint.
 Also, 
$$\sfT(b)=exp_u(\sfT(D_u(b_1,\ldots,b_j)),\sfS_{j+1},\ldots,\sfS_n)$$
 for some $\sfS_{j+1},\ldots,\sfS_n\in M2T$ such that $\sfS_k\cong \sfT(b_k)$ for $j+1\leq k \leq n$ and $\sfT(D_u(b_1,\ldots,b_n)),\sfS_{j+1},\ldots,\sfS_n$ are pairwise disjoint.
 
We will define $\sf A$ to be $\sfT(D_u(a_1,\ldots,a_i))$ in case $i\not=0$ and define $\sf B$ to be $\sfT(D_u(b_1,\ldots,b_j))$ in case $j\not=0$.
 
 The assumption that $\sfT(a)\, \preq^c\ \sfT(b)$ implies 
 $$\oplus({\sf A},\sfT_{i+1},\ldots,\sfT_m) \, \preq^c \ \oplus({\sf B},\sfS_{j+1},\ldots,\sfS_n)$$
 The rest of the argument for this case is similar to that for Subcase 1.3 and is omitted.
 \qed

\blankline

We will need the following operation transforming sequences of double forests into a double tree by adding an element which becomes a simultaneous root. 

\begin{dfn}
Assume ${\bf P}_1,\ldots,{\bf P}_m$ are pairwise disjoint double forests. 
Define $\rho(\bfP_1,\ldots,\bfP_m)$ to be a structure $(X,\leqone,\leqtwo)$ which satisfies the following conditions.
\begin{enumerate}
\item
$X=|\bfP_1|\cup \cdots \cup |\bfP_n|\cup \{r\}$ where $r\not\in |\bfP_i|$ for $i=0,\ldots,n$.
\item
Assume $1\leq i\leq n$ and $x\in |\bfP_i|$.
For $k=1,2$ and all $y\in X$, $x\leq_k y$ iff $x\leq_k^{{\bf P}_i}y$.
\item
For $k=1,2$ and all $y\in X$, $r\leq_k y$.
\end{enumerate}
 If $\bfP_1,\ldots,\bfP_m$ are double forests which are not pairwise disjoint then define $\rho(\bfP_1,\ldots,\bfP_m)$ to be $\rho(\bfQ_1,\ldots,\bfQ_m)$ where $\bfQ_1,\ldots,\bfQ_m$ are chosen to be pairwise disjoint and so that $\bfQ_i\cong \bfP_i$ for $i=1,\ldots,m$.
\end{dfn}

\begin{lem}
\label{ldoubletree}
If $\bfP_1,\ldots,\bfP_m$ are double forests then $\rho(\bfP_1,\ldots,\bfP_m)$ is a double tree.
\end{lem}
{\bf Proof.} 
Clear.
\qed

\begin{cor}
\label{cmain}
Assume $a=(a_1,\ldots,a_m)$ and $b=(b_1,\ldots,b_n)$ are in $OT(\omega)$ and have order $0$.
If $\rho(\sfT(a_1),\ldots,\sfT(a_m))\, \preq^c \ \rho(\sfT(b_1),\ldots,\sfT(b_n))$ then $a\preq b$.
\end{cor}
{\bf Proof.}
Recall that we have identified elements of $M2F$ of order 0 with double forests.

Assume $\rho(\sfT(a_1),\ldots,\sfT(a_m))\, \preq^c \ \rho(\sfT(b_1),\ldots,\sfT(b_n))$.

There are pairwise disjoint $\sfT_1,\ldots,\sfT_m$ with $\sfT_i\cong \sfT(a_i)$ for $i=1,\ldots,m$ such that $\rho(\sfT(a_1),\ldots,\sfT(a_m))=\rho(\sfT_1,\ldots,\sfT_m)$.
There are pairwise disjoint $\sfS_1,\ldots,\sfS_n$ with $\sfS_j\cong \sfT(b_j)$ for $j=1,\ldots,n$ such that $\rho(\sfT(b_1),\ldots,\sfT(b_n))=\rho(\sfS_1,\ldots,\sfS_n)$.
Let $h$ be a covering of $\rho(\sfT_1,\ldots,\sfT_m)$ into $\rho(\sfS_1,\ldots,\sfS_n)$.
The restriction of $h$ is clearly a covering of  $\oplus(\sfT_1,\ldots,\sfT_m)$ into $\oplus(\sfS_1,\ldots,\sfS_n)$.
We can now argue as in Subcase 1.3 of Theorem \ref{tmain}.

For $1\leq i \leq m$, there exists $p(i)$ such that $h$ maps $|\sfT_i|$ into $|\sfS_{p(i)}|$.
By Theorem \ref{tmain}, $a_i\preq b_{p(i)}$.
Moreover, if $1\leq i,j \leq n$, $i\not=j$ and $p(i)=p(j)$ then $h$ cannot map $|\sfT_i|$ onto $|\sfS_{p(i)}|$ implying $a_i\not=b_{p(j)}$.
By the observation preceding Theorem \ref{tmain}, $a\preq b$.
\qed

\section{From WQO to WO}

The previous sections provide the means to prove the following theorem along with several variations. 
Recall that the proofs in the previous two sections can be carried out in ${\bf RCA}_0$.
Except where specified otherwise, the arguments in this section and the next can also be formalized in ${\bf RCA}_0$.

\begin{thm}{\bf (RCA$_0$)}
\label{treduction1}
If the collection of double trees is wqo under covering then $\{a\in OT \, : \, a \prec D_0D_\omega0\}$ is well-ordered by $\prec$.
\end{thm}

The proof of the theorem will use the following lemma.

\begin{lem}
\label{lD0Domega}
Assume $a\in OT$.
\begin{enumerate}
\item
If $D_0a\in OT$ and $ord(a)\leq n<\omega$ then $D_0a\in OT(n+1)$. 
\item
If $n$ is a natural number and $a \prec D_0D_{n+1}0$ then $ord(a)=0$ and $a\in OT(n+1)$.
\item
If $a\prec D_0D_\omega 0$ then $ord(a)=0$ and $a\in OT(\omega)$.
\end{enumerate}
\end{lem}
{\bf Proof.}
Notice that $b\prec D_{n+1}0$ iff $ord(b)\leq n$ for all $b\in OT$.
A straightforward induction shows that for all $b\in OT$, if $ord(b)\leq n$ and each element of $G_0b$ has order at most $n$ then $b\in OT(n+1)$.

For part 1, assume $D_0a\in OT$ and $ord(a)\leq n < \omega$.
Since $D_0a\in OT$, $G_0a\prec a$ implying each element of $G_0a$ has order at most $n$.
By the observation in the previous paragraph, $a\in OT(n+1)$ implying $D_0a\in OT(n+1)$.

For part 2,  suppose $a\prec D_0D_{n+1}0$ where $n<\omega$.
Since $0\in OT(n+1)$, we may assume $a\not=0$.
We must have $ord(a)=0$. 
If $a=(a_1,\ldots,a_m)$ where $1<m$, it suffices to show $a_i\in OT(n+1)$ for $i=1,\ldots,m$.
Therefore, we may assume $a$ is a principal term. 
Since $ord(a)=0$, $a=D_0b$ for some $b$.
Since $D_0b=a\prec D_0 D_{n+1}0$, $b\prec D_{n+1}0$ implying $ord(b)\leq n$.
By part 1, $a=D_0b \in OT(n+1)$.

A similar argument establishes part 3.
\qed

\blankline

\noindent{\bf Proof of Theorem \ref{treduction1}.}
Assume the collection of double trees is wqo under covering.
Suppose $a_0,a_1,\ldots,a_i,\ldots$ is an infinite sequence of elements of $\{a\in OT \, : \, a\prec D_0D_\omega0\}$.
We must find $i<j$ such that $a_i\preq a_j$.
By part 3 of Lemma \ref{lD0Domega}, $ord(a_i)=0$ and $a_i\in OT(\omega)$ for all $i<\omega$.
We may assume $a_i\not=0$ for all $i<\omega$ (if $a_i=0$ then $a_i\preq a_{i+1}$).
For $i<\omega$, define $\sfT_i$ as follows.
There are principal $b_1,\ldots,b_k$ such that $a_i=(b_1,\ldots,b_k)$ (allowing the possibility $k=1$).
Since $ord(a_i)=0$, we must have $ord(b_j)=0$ for $j=1,\ldots,k$.
Hence, $\sfT(b_j)$ has order $0$ for $j=1,\ldots,k$.
Let $\sfT_i=\rho(\sfT(b_1),\ldots,\sfT(b_k))$.
By Lemma \ref{ldoubletree}, each $\sfT_i$ is a double tree.
By assumption, there are $i<j$ such that $\sfT_i\, \preq^c  \ \sfT_j$.
By  Corollary \ref{cmain}, $a_i\preq a_j$.
\qed

\blankline

Assume ${\cal Q} = (Q,\leq)$ is a quasiordering. 
We will write $WQO({\cal Q})$ to indicate that $\cal Q$ is a wqo.
We will write $PRWQO({\cal Q})$ to indicate there are no primitive recursive bad sequences in $\cal Q$.

When $\cal Q$ is a linear ordering, $\cal Q$ is a well-ordering iff $WQO({\cal Q})$ and $\cal Q$ has no primitive recursive descending sequences iff $PRWQO({\cal Q})$. In this case, we write  $WO({\cal Q})$ for $WQO({\cal Q})$ and $PRWO({\cal Q})$ for $PRWQO({\cal Q})$. 

When $a\in OT$, $X$ is the set of $b\in OT$ such that $b\prec a$ and $\alpha=o(a)$ (where the operation $a\mapsto o(a)$ is defined in [\ref{rbu}]), we will write $WO(\alpha)$ and $PRWO(\alpha)$  for $WO((X,\preq))$ and $PRWO((X,\preq))$ respectively. 
The following calculations are from [\ref{rbu}]: $o(D_0D_\omega0)=\psi_0\Omega_\omega$, $o(D_0D_{n+1}0)=\psi_0\varepsilon_{\Omega_n+1}$ for $1\leq n<\omega$ and $o(D_0D_10)=\varepsilon_0$.

The reader unfamiliar with $[\ref{rbu}]$ may simply view the notations $WO(\alpha)$ and $PRWO(\alpha)$ from the previous paragraph as abbreviations.

We will write $DTC$ for the partial ordering of double trees under covering.
For $n\in\omega$, $DTC(n)$ is the restriction of $DTC$ to the collection of double trees of height at most $n$.
We write $TC$ for partial ordering of finite trees under covering (where coverings between trees are defined analagously to coverings between double forests i.e. a covering preserves order upward).

\begin{lem}
\label{ltc}
$DTC(1)$ is isomorphic to $TC$.
\end{lem}
{\bf Proof.}
Notice that for any $(X,\leq_1,\leq_2)\in DT(1)$, $x\leq_2 y$ iff $x$ is the root of $(X,\leq_1,\leq_2)$ for all $x,y\in X$.
This implies that the map $(X,\leq_1,\leq_2)\mapsto (X,\leq_1)$ from $DTC(1)$ to $TC$ is a bijection.

Suppose ${\bf P}_1=(X_1,\leq^1_1,\leq^1_2)$ and ${\bf P}_2=(X_2,\leq^2_1,\leq^2_2)$ are double trees with roots $r_1$ and $r_2$ respectively. 
Clearly, if $h:X_1\rightarrow X_2$ is a covering of ${\bf P}_1$ into ${\bf P}_2$ then it is also a covering of $(X_1,\leq^1_1)$ into $(X_2,\leq^2_1)$.
Now suppose $h:X_1\rightarrow X_2$ is a covering of $(X_1,\leq^1_1)$ into $(X_2,\leq^2_1)$.
A simple argument shows that if we modify $h$ by mapping $r_1$ to $r_2$ the result is a covering of ${\bf P}_1$ into ${\bf P}_2$.
\qed

\begin{thm}{\bf (RCA$_0$)}
\label{treduction2}
\begin{enumerate}
\item
$WQO(DTC) \ \Longrightarrow \ WO(\psi_0\Omega_\omega)$
\item
For $1\leq n <\omega$, $WQO(DTC(n+1)) \ \Longrightarrow \ WO(\psi_0\varepsilon_{\Omega_n+1})$
\item
$WQO(TC) \ \Longrightarrow \ WO(\varepsilon_0)$
\item
$PRWQO(DTC) \ \Longrightarrow \ PRWO(\psi_0\Omega_\omega)$
\item
For $1\leq n <\omega$, $PRWQO(DTC(n+1)) \ \Longrightarrow \ PRWO(\psi_0\varepsilon_{\Omega_n+1})$
\item
$PRWQO(TC) \ \Longrightarrow \ PRWO(\varepsilon_0)$
\end{enumerate}
\end{thm}

We remark that the proofs for parts 3 and 6 are fairly direct and require little of the development of the previous sections. 

The functions defined on $OT$, $M2T$ or $M2F$ in the previous two sections are of low complexity (after making natural choices to be explicit regarding outputs when necessary). 
This is clear from the definitions along with the descriptive lemmas.
For our present purposes, we  only need to notice they are primitive recursive.

\begin{lem}
If $0\not= a\in OT(\omega)$, $ord(a)=0$  and $n$ is maximal such that $D_n$ occurs in $a$ then $\rho(\sfT(a_1),\ldots,\sfT(a_k))$ is a double tree of height $n+1$
\end{lem}
{\bf Proof.}
For $\sfT\in M2F$, define the {\it height} of \sfT\ to be the largest number of the form $v+k$ where there is a chain in $\leq_2^{\sf T}$ of size $k+1$ all of whose elements have order $v$. 
Notice that if $\sfT$ has order 0 then this definition of height agrees with our previous definition: the largest $m$ such there is a chain in $\leq_2^{\sf T}$ of size $m+1$.

\halfblankline

{\bf Claim.} If $\sfT_1,\ldots,\sfT_m \in M2T$ is order descending and $u<\omega$ then the height of $\Psi_u(\sfT_1,\ldots,\sfT_m)$ is the maximum of $u$ and the heights of $\sfT_i$ for $i=1,\ldots,n$.

Straightforward induction on the cardinality of $\Psi_u(\sfT_1,\ldots,\sfT_m)$ noticing that $coll$ may increase the sizes of chains but reduces the corresponding order and $exp_u$ does not increase the sizes of chains and  preserves order while adding a new root of order $u$.

\halfblankline

{\bf Claim 2.} For all principal $a\in OT(\omega)$, if $n$ is maximal such that $D_n$ occurs in $a$ then the height of $\sfT(a)$ is $n$.

The claim follows from the previous claim by a straightforward induction noting that the height of $\sfT(D_u0)$ is $u$ for $u<\omega$.

\halfblankline

The lemma follows immediately from Claim 2 noting that $\rho$ increases height by 1.
\qed

\blankline

\noindent{\bf Proof of Theorem \ref{treduction2}.} Part 1 is simply a restatement of Theorem \ref{treduction1} using abbreviated notation. The proofs of parts 2 through 6 follow the same lines with slight modifications.

For part 2, modify the proof of Theorem \ref{treduction1} by using part 2 of Lemma \ref{lD0Domega} rather than part 3 and using the previous lemma to see that the double trees $\sfT_i$ have height at most $n+1$.

For part 3, notice the proof of part 2 also works for $n=0$ to show $WQO(DTC(1))$ implies that $\{a\in OT \, : \, a\prec D_0D_10\}$ is well-ordered i.e. $WO(\varepsilon_0)$.
Part 3 now follows from the fact $DTC(1)$ is isomorphic to $TC$.

To prove parts 4 through 6, modify the proofs of  parts 1 through 3 respectively by noticing that if $a_0,a_1,\ldots,a_i,\ldots$ is a primitive recursive sequence then so is $\sfT_0,\sfT_1,\ldots,\sfT_i,\ldots$.
\qed

\blankline

Assume ${\cal Q} = (Q,\leq^{\cal Q})$ is a quasiordering.
Given a {\it norm} $q\mapsto ||q||$ from $Q$ into $\omega$, we will write $LWQO({\cal Q})$, (as in [\ref{rsi85}]), to indicate there is no bad sequence $q_0,q_1,\ldots,q_i,\ldots$ such that the sequence $||q_0||,||q_1||,\ldots,||q_i||,\ldots$ of natural numbers is bounded by a linear function.
As observed in [\ref{rsi85}], when  $\{q\in Q \, : \, ||q||\leq n\}$ is finite for all $n\in \omega$, $LWQO({\cal Q})$ is equivalent, in ${\bf WKL}_0$, to the following $\Pi^0_2$ statement: $\forall c\in \omega \ \exists k\in \omega$ such that if $q_0,\ldots,q_k$ are elements of $Q$ with $||q_i||\leq c\cdot(i+1)$ for $i=0,\ldots,k$ then there exist $i<j$ such that $q_i\leq^{\cal Q}q_j$.
We will use this latter statement as our official definition of $LWQO({\cal Q})$ in ${\bf RCA}_0$.

When $\cal Q$ is a linear ordering, we will write  $LWO(\cal Q)$ for $LWQO(\cal Q)$.

We will fix norms on $M2F$ and $OT$ such that $||\sfT||$ is the cardinality of $|\sfT|$ for $\sfT\in M2F$ and $||a||$ is the length of $a$ for $a\in OT$.

\begin{thm}\
\label{treduction3}
\begin{enumerate}
\item
$({\bf RCA}_0)$ For $1\leq n < \omega$, $LWQO(DTC(n+1)) \ \Longrightarrow \ LWO(\psi_0\varepsilon_{\Omega_n+1})$
\item
$({\bf RCA}_0)$ $LWQO(TC) \ \Longrightarrow \ LWO(\varepsilon_0)$
\item
$({\bf ACA}_0)$ For $1\leq n < \omega$, $LWQO(DTC(n+1)) \ \Longrightarrow \ PRWO(\psi_0\varepsilon_{\Omega_n+1})$
\item
$({\bf ACA}_0)$ $LWQO(TC) \ \Longrightarrow \ PRWO(\varepsilon_0)$
\item
$({\bf ACA}_0)$ $LWQO(DTC) \ \Longrightarrow \ PRWO(\psi_0\Omega_\omega)$
\end{enumerate}
\end{thm}

We remark that $LWO(\psi_0\Omega_\omega)$ is false.
Consider the descending sequences $D_0D_{n+1}0,D_0D_n0,\ldots,D_0D_10$ where all the terms have length 3.

The proof of the theorem will use the following lemma.

\begin{lem}
Assume $n<\omega$. 
If $a\in OT(n+1)$ is principal then $||\sfT(a)||\leq (n+1)\cdot ||a||$. 
\end{lem}
{\bf Proof.}
Compare the following claim to Lemma \ref{lPsipicture}.

\halfblankline 

{\bf Claim.} If $\sfT_1,\ldots,\sfT_k\in M2T$ is order descending where $u<ord(\sfT_k)$ then 
$$||\Psi_u(\sfT_1,\ldots,\sfT_k)||= ||\sfT_1||+\cdots+||\sfT_k||+(ord(\sfT_1)-u)$$

The proof of the claim is a straightforward induction on $ord(\sfT_1)-u$ using Definition \ref{dpsi1}.

\halfblankline

The claim implies that if $\sfT_1,\ldots,\sfT_k\in M2T$ is order descending and $ord(\sfT_1)\leq n$ then 
$$||\Psi_u(\sfT_1,\ldots,\sfT_k)||\leq ||\sfT_1||+\cdots+||\sfT_k||+(n+1)$$
This follows immediately from the claim if $u<ord(\sfT_k)$. Otherwise, use the claim and refer to Definition \ref{dpsi2}.

The lemma now follows by induction on $a\in OT(n+1)$.
\qed

\blankline

\noindent{\bf Proof of Theorem \ref{treduction3}.}
For part 1, assume $1\leq n<\omega$. 
 Once again, the proof is similar to the proof of Theorem \ref{treduction1}.
 
Suppose $c>0$. 
There exists $k<\omega$ such that for all double trees $\sfT_0,\ldots,\sfT_k$ of height at most $n+1$ with $||\sfT_i||\leq (n+2)\cdot c \cdot (i+1)$ for $i=0,\ldots,k$, there exists $i<j$ such that $\sfT_i\preq^c \sfT_j$.
Assume $a_i\prec D_0D_{n+1}0$ with $||a_i||\leq c\cdot(i+1)$ for $i=0,\ldots,k+1$.
By part 2 of Lemma \ref{lD0Domega}, $a_i\in OT(n+1)$ and $ord(a_i)=0$ for $i=0,\ldots,k+1$. 
We will show there exist $i<j$ such that $a_i\preq a_j$.
We may assume $a_i\not = 0$ for $i\leq k$.

For $i=0,\ldots,k$, define $\sfT_i$ as follows.
There exist principal $b_1,\ldots,b_m\in OT(n+1)$ such that $a_i=(b_1,\ldots,b_m)$ (allowing the possibility $m=1$). 
Let $\sfT_i=\rho(\sfT(b_1),\ldots,\sfT(b_m))$.
Using the previous lemma, for $i\leq k$
\begin{tabbing}
\hspace{1in}  $||\sfT_i||$ \= = \ $||\sfT(b_1)||+\cdots+||\sfT(b_m)||+1$ \\
                              \> $\leq \ (n+1)\cdot ||b_1||+\cdots + (n+1)\cdot||b_m)||+1$ \\
                              \> = \ $(n+1)\cdot(||b_1||+\cdots +||b_m||)+1$ \\
                              \> $\leq \ (n+1)\cdot||a_i||+1$ \\
                              \> $\leq \ (n+2)\cdot||a_i||$ \\
                              \> $\leq \ (n+2)\cdot c \cdot (i+1)$
\end{tabbing}
By choice of $k$, there exist $i<j$ such that $\sfT_i\preq^c \sfT_j$.
By Corollary \ref{cmain}, $a_i\preq a_j$.

The argument for part 1 also works for $n=0$ to show $LWQO(DT(1))$ implies $LWO(\varepsilon)$. 
Since $TC\cong DTC(1)$ by Lemma \ref{ltc}, part 2 follows.

Parts 3 and 4 follow from parts 1 and 2 using Theorem 3.5 from [\ref{rsi85}], having assumption ${\bf ACA}_0$,  which has as special cases that $LWO(\psi_0\varepsilon_{\Omega_n+1})$ implies $PRWO(\psi_0\varepsilon_{\Omega_n+1})$ for $1\leq n<\omega$ and $LWO(\varepsilon_0)$ implies $PRWO(\varepsilon_0)$.

For part 5, assume $LWQO(DTC)$.
By part 3, $PRWO(\psi_0\varepsilon_{\Omega_n+1})$ whenever $1\leq n < \omega$.
Since $D_0D_{n+1}0$ ($n\in\omega$) is cofinal in $D_0D_\omega0$ by part 3 of Lemma \ref{lD0Domega}, this implies $PRWO(\psi\Omega_\omega)$.
\qed

\section{Applications: Strength and Independence}

For a theory $\bf T$ and a collection of formulas $\Phi$, the {\bf uniform $\Phi$ reflection principle for $\bf T$} is the collection of formulas formalizing the statements 
\begin{center}
{\it For every natural number $n$, if $\varphi(n)$ is provable in $\bf T$ then $\varphi(n)$ is true.} 
\end{center}
where $\varphi$ is a formula in $\Phi$ with at most one free variable.

Generally, a proof-theoretic analysis of a theory $\bf T$ showing that the ordinal of $\bf T$ is $\alpha$ also shows, though possibly not stated explicitly,  that the following are provable in ${\bf ACA}_0$:
\begin{itemize}
\item[$(*)$] $WO(\alpha)$ is equivalent to the uniform $\Pi^1_1$ reflection principal for $\bf T$.
\item[$(**)$] $PRWO(\alpha)$ is equivalent to the uniform $\Pi^0_2$ reflection principal for $\bf T$.
\end{itemize}
In particular, $(*)$ and $(**)$ hold for $\Pi^1_1-{\bf CA}_0$ and  $\psi_0\Omega_\omega$ as well as ${\bf ID}_n$ and $\psi_0\varepsilon_{\Omega_n+1}$ when $1\leq n <\omega$ (see [\ref{rratowe}]).
 Rathjen and Weiermann [\ref{rrawe}] gave an ordinal analysis providing an appropriate instance of $(*)$ to calibrate the strength of Kruskal's Theorem.
 
 \begin{thm} \
 \label{tstrength}
 \begin{enumerate}
 \item $({\bf ACA}_0)$ 
$WQO(DTC)$ is equivalent to the uniform $\Pi^1_1$ reflection principal for $\Pi^1_1-{\bf CA}_0$.
 \item $({\bf ACA}_0)$ The following are equivalent.
 \begin{enumerate}
 \item $PRWQO(DTC)$
 \item $LWQO(DTC)$
 \item The uniform $\Pi^0_2$ reflection principal for $\Pi^1_1-{\bf CA}_0$.
 \end{enumerate}
 \end{enumerate}
 \end{thm}
 
The final section of [\ref{rca}] establishes that the Double Kruskal Theorem follows from ${\bf RCA}_0$ with the additional assumption of the uniform $\Pi^1_1$ reflection principal for ${\bf KP}\ell_0$.
The proof is based on showing that  
 \begin{itemize}
 \item[]
 $(\dagger)$ \ \ 
 For all $n\in\omega$, $WQO(DTC(n))$ is provable in ${\bf KP}\ell_0$.
 \end{itemize}
follows from ${\bf RCA}_0$.
The proof of Theorem \ref{tstrength} will use a variant which says
 \begin{itemize}
 \item[]
 $(***)$ \ \ 
 For all $n\in\omega$, $WQO(DTC(n))$ is provable in $\Pi^1_1-{\bf CA}_0$.
 \end{itemize}
 follows from ${\bf RCA}_0$.
 This is not surprising since ${\bf KP}\ell_0$ is a conservative extension of $\Pi^1_1-{\bf CA}_0$ (e.g. see Chapter 7 of [\ref{rsi10}], especially Exercise VII.3.36).
 In fact, the proof that $(\dagger)$ follows from ${\bf RCA}_0$ also shows, with only cosmetic changes, that $(***)$ follows from ${\bf RCA}_0$.
 
 We will also need the following observations.
 
 \begin{lem}\
 \begin{enumerate}
 \item
$WQO(DTC) \ \Longleftrightarrow \ \forall n \, WQO(DTC(n))$
 \item
$PRWQO(DTC) \ \Longleftrightarrow \ \forall n \, PRWQO(DTC(n))$
 \item
$LWQO(DTC) \ \Longleftrightarrow \ \forall n \, LWQO(DTC(n))$
 \end{enumerate}
 \end{lem}
{\bf Proof.}
The forward direction of each part is obvious.
The reverse directions use the following claim.

\halfblankline

{\bf Claim.} Assume ${\bf P}$ and ${\bf Q}$ are double trees.
If the height of $\bf Q$ is at least $||{\bf P}||-1$  then ${\bf P}\preq^c \ {\bf Q}$.

 Assume the height of $\bf Q$ is at least $||{\bf P}||-1$.
 Suppose ${\bf P}=(X,\leq_1,\leq_2)$ and ${\bf Q}=(Y,\preq_1,\preq_2)$.
 Let $C$ be a chain in $(Y,\preq_2)$ of size $||{\bf P}||$ and extend $\leq_1$ to a linear ordering $\leq_1'$ of $X$.
 Let $h$ map $X$ into $C$ so as to preserve order between $\leq_1'$ and $\preq_2$.
Clearly, $h$ is a covering of $\bf P$ into $\bf Q$.
 
 \halfblankline

 To prove the reverse direction of part 1, assume $WQO(DTC(n))$ for all $n\in \omega$ and let ${\bf P}_0,\ldots,{\bf P}_i,\ldots$ be an infinite sequence of double trees. 
 We will show there are $i<j$ such that ${\bf P}_i\preq^c \ {\bf P}_j$.
By the claim, we may assume the height of ${\bf P}_i$ is less than $||{\bf P}_0||-1$ for all $i\geq 1$. 
By $WQO(DTC(n))$ where  $n=||{\bf P}_0||-2$, there must be $i,j\in\omega$ with $1\leq i <j$ such that ${\bf P}_i\preq^c\ {\bf P}_j$.

The proof of the reverse direction of part 2 is identical except we assume $PRWQO(DTC(n))$ for all $n$ and ${\bf P}_0,\ldots,{\bf P}_i,\ldots$ is a primitive recursive sequence of double trees.

To prove the reverse direction of part 3, assume $LWQO(DTC(n))$ for all $n\in\omega$.
Suppose $c>0$. 
There exists $k\in\omega$ such that for any sequence ${\bf Q}_0,\ldots,{\bf Q}_k$ of double trees of height less than $c-1$ with $||{\bf Q}_i||\leq 2c(i+1)$ for $i=0,\ldots,k$ there exist $i<j$ such that ${\bf Q}_i\preq^c \ {\bf Q}_j$.
Suppose ${\bf P}_0,\ldots,{\bf P}_{k+1}$ is a sequence of double trees with $||{\bf P}_i||\leq c(i+1)$ for $i=0,\ldots,k+1$. 
By assumption, $||{\bf P}_0||\leq c$.
By the claim, we may assume the height of ${\bf P}_{i+1}$ is less than $||{\bf P}_0||-1$ and, hence, less than $c-1$ for $i=0,\ldots,k$.
Consider the sequence ${\bf P}_{i+1}$ ($i=0,\ldots,k$) and notice $||{\bf P}_{i+1}||\leq c(i+2)\leq 2c(i+1)$ for $i=0,\ldots,k$.
By the choice of $k$, there exist $i<j$ such that ${\bf P}_{i+1}\preq^c \ {\bf P}_{j+1}$. 
\qed

\blankline

 \noindent{\bf Proof of Theorem \ref{tstrength}.}
 We use the instances of $(*)$ and $(**)$ with ${\bf T}=\Pi^1_1-{\bf CA}_0$ and $\alpha=\psi_0\Omega_\omega$.
 
 The forward direction of part 1 follows from the fact that $WQO(DTC)$ implies $WO(\psi_0\Omega_\omega)$ (part 1 of Theorem \ref{treduction2}) and $(*)$.
 
 For the reverse direction of part 1, we assume the uniform $\Pi^1_1$ reflection principal for $\Pi^1_1-{\bf CA}_0$. By $(***)$, $\Pi^1_1-{\bf CA}_0\vdash WQO(DTC(n))$ for each $n\in \omega$.
 Since $WQO(DTC(n))$ is clearly equivalent to a $\Pi^1_1$ statement, this implies $WQO(DTC(n))$ for all $n\in\omega$.
 By part 1 of the previous lemma, we have $WQO(DTC)$.

  For part 2, we will show $(a) \Leftrightarrow (c)$ and $(b)\Leftrightarrow (c)$.
 
 The implication $(a)\Rightarrow (c)$ follows from  the fact that $PRWQO(DTC)$ implies $PRWO(\psi_0\Omega_\omega)$ (part 4 of Theorem \ref{treduction2}) and $(**)$.
 
 The implication $(b)\Rightarrow (c)$ follows from  the fact that $LWQO(DTC)$ implies $PRWO(\psi_0\Omega_\omega)$ (part 5 of Theorem \ref{treduction3}) and $(**)$.
 
 The proofs that $(c)$ implies $(a)$ and $(b)$  are similar to the proof of the reverse direction of part 1.
  Assume the uniform $\Pi^0_2$ reflection principal for $\Pi^1_1-{\bf CA}_0$. By $(***)$, $\Pi^1_1-{\bf CA}_0\vdash WQO(DTC(n))$ for each $n\in \omega$.
  
 To prove $(a)$, notice that $\Pi^1_1-{\bf CA}_0\vdash PRWQO(DTC(n))$ for $n\in\omega$ 
 (since it is provable in ${\bf RCA}_0$ and, hence, in $\Pi^1_1-{\bf CA}_0$ that $WQO(DTC(n))$ implies $PRWQO(DTC(n))$).
 Since $PRWQO(DTC(n))$ is equivalent to a $\Pi^0_2$ statement, this implies $PRWQO(DTC(n))$ holds for all $n\in\omega$.
 By part 2 of the previous lemma, we have $PRWQO(DTC)$.
 
The proof of $(b)$ is similar.
Notice that $\Pi^1_1-{\bf CA}_0\vdash LWQO(DTC(n))$ for $n\in\omega$ 
 (since it is provable in ${\bf WKL}_0$ and, hence, in $\Pi^1_1-{\bf CA}_0$ that $WQO(DTC(n))$ implies $LWQO(DTC(n))$).
 Since $LWQO(DTC(n))$ is equivalent to a $\Pi^0_2$ statement, this implies $LWQO(DTC(n))$ holds for all $n\in\omega$.
 By part 3 of the previous lemma, we have $LWQO(DTC)$.

\qed

\blankline

The previous theorem also holds when $DTC$ is replaced by $TC$ and $\Pi^1_1-{\bf CA}_0$ is replaced by ${\bf ACA}_0$. One can use $(*)$ and $(**)$ with $\alpha=\varepsilon_0$ and replace $DTC(n)$ by the collection of trees of height at most $n$ in the proof above, but a direct proof avoids many of the difficulties in earlier sections.
We remark that the proposition that $TC$ is wqo under covering follows from the bounded version of Kruskal's Theorem.

\begin{cor}
$({\bf RCA}_0)$ If $\Pi^1_1-{\bf CA}_0$ is consistent then neither $PRWQO(DTC)$ nor $LWQO(DTC)$ is provable in $\Pi^1_1-{\bf CA}_0$ and, hence, $WQO(DTC)$ is not provable in $\Pi^1_1-{\bf CA}_0$.
\end{cor}
{\bf Proof.}
By part 2 of the theorem and G\"odel's Incompleteness Theorem.
\qed

 \begin{thm}
 \label{tstrength2}
 Assume $1\leq n <\omega$.
 \begin{enumerate}
 \item
  $({\bf ACA}_0)$
$PRWQO(DTC(n+1))$ implies the uniform $\Pi^0_2$ reflection principal for ${\bf ID}_n$.
 \item 
  $({\bf ACA}_0)$
 $LWQO(DTC(n+1))$ implies the uniform $\Pi^0_2$ reflection principal for ${\bf ID}_n$.
 \end{enumerate}
 \end{thm}
 {\bf Proof.}
We use the instance of $(**)$ with ${\bf T}={\bf ID}_n$ and $\alpha=\psi_0\varepsilon_{\Omega_n+1}$.

Part 1 follows from the fact that $PRWQO(DTC(n+1))$ implies $PRWO(\psi_0\varepsilon_{\Omega_n+1})$ (part 5 of Theorem \ref{treduction2}) and $(**)$.

Part 2 follows from the fact that $LWQO(DTC(n+1))$ implies $PRWO(\psi_0\varepsilon_{\Omega_n+1})$ (part 3 of Theorem \ref{treduction3}) and $(**)$.
\qed

\blankline

We expect that the converses of the parts of the previous theorem also hold establishing the analogue of part (b) of Theorem \ref{tstrength}. 

\begin{cor}
Assume $1\leq n <\omega$.  \\
$({\bf RCA}_0)$ If ${\bf ID}_n$ is consistent then neither $PRWQO(DTC(n+1))$ nor $LWQO(DTC(n+1))$ is provable in ${\bf ID}_n$. 
\end{cor}
{\bf Proof.}
By the theorem and G\"odel's Incompleteness Theorem.
\qed

\blankline

The results of this section concerning $LWQO({\bf Q})$ when $\bf Q$ is one of $TC$, $DTC(n+1)$ or $DTC$ open the door to questions in the area of phase transition as developed by A. Weiermann.

%: References

\begin{center}
REFERENCES
\end{center}

\begin{enumerate}
\item
\label{rbu}
W. Buchholz, {\it A new system of proof-theoretic ordinal functions,} Annals of Pure and Applied Logic 32 (1986), pp. 195-207.
\item
\label{rca}
T. Carlson, {\it Generalizing Kruskal's Theorem to pairs of cohabitating trees,}  Archive for Mathematical Logic 55 (2016), pp. 37-48.
\item
\label{rja}
G. J\"ager, {\bf Theories for Admissible Sets: A Unifying Approach to Proof Theory,} Bibliopolis, Napoli, 1986, v + 167 pages.
\item
\label{rkr}
J. B.  Kruskal, {\it Well-quasi-ordering, the tree theorem, and Vazsonyi's conjecture}, Transactions of the American Mathematical Society 95 (1960), pp. 210-225. 
\item
\label{rla76}
R. Laver, {\it Well-quasi-orderings and sets of finite sequences}, Mathematical Proceedings of the Cambridge Philosophical Society 79 (1976), pp. 1-10.
\item
\label{rratowe}
M. Rathjen, M. Toppel and A. Weiermann, {\it Ordinal analysis, proof-theoretic reductions and conservativity,} in preparation.
\item
\label{rrawe}
M. Rathjen and A. Weiermann, {\it Proof-theoretic investigations on Kruskal's theorem,} Annals of Pure and Applied Logic 60 (1993), pp. 49-88.
\item
\label{rsi85}
S. G. Simpson, {\it Nonprovability of certain combinatorial properties of finite trees,} {\bf Harvey Friedman's Research on the Foundations of Mathematics,} (L. A. Harrington, M. D. Morley, A. Scedrov, and S. G. Simpson, editors), North-Holland, Amsterdam, 1985, pp. 87-117.
\item
\label{rsi10}
S. G. Simpson, {\bf Subsystems of Second Order Arithmetic, 2nd Edition,} Perspectives in Logic, Cambridge University Press, Cambridge, GBR, 2010, xvi + 444 pages. 
\end{enumerate}

\end{document}